# Exact expressions for the weights used in least-squares regression estimation for the log-logistic and Weibull distribution

J. Martin van Zyl

*Department of Mathematical Statistics and Actuarial Science, University of the Free State, PO Box 339, Bloemfontein, South Africa*

Abstract: Estimation for the log-logistic and Weibull distributions can be performed by using the equations used for probability plotting. The equations leads to highly heteroscedastic regression. Exact expressions for the variances of the residuals are derived which can be used to perform weighted regression. In large samples maximum likelihood performs best, but it is shown that in smaller samples the weighted regression outperforms maximum likelihood estimation with respect to bias and mean square error.

**1 Introduction**

The Weibull and log-logistic distributions are of importance in survival analysis where both distributions have hazard functions which have simple forms and can have several shapes. One of the methods of estimation is based on probability plotting. This leads to a linear regression involving the order statistics which have different variances and heteroscedasticity is a problem. In this work the distributions of the residuals and exact expressions are derived for the variances in the regression which makes generalized least squares (GLS) regression possible. The residuals variances are independent of the underlying distribution in both cases, and can be used for any parameters.

For the log-logistic regression the residuals are logistic distributed and in the case of the Weibull regression the distribution of the residuals is a Gumbel distribution for the minimum. In both cases the full covariance matrices have determinants which are for all practical purposes zero, even though the covariance matrixes are of full rank. Since the determinant is extremely small, that is less than say $10^{-5}$ even in small samples, weighted regression will be used, involving a diagonal covariance matrix for GLS. It is also shown that by using a very simple simulation involving uniform random variables, the covariance matrix can be calculated with a high accuracy.

It was found that the weighted regression outperforms maximum likelihood (ML) estimation in small samples, and especially with respect to bias. The other method of estimation which is frequently used for estimation of the parameters of these two distributions is the probability-weighted (PWM) method of estimation introduced by Greenwood *et al.* (1979). Several studies were conducted to evaluate the performance of this method and it will not be considered in this work. Some references involving PWM estimation are the papers by Hosking (1990), Askhar and Mahdi (2003), Shoukri *et al.* (1988), Whalen *et al.* (2004), Genschel and Meeker (2010). PWM estimation performs well where data are not heavy-tailed and in large samples ML estimation performs better that this method.

Van Zyl and Schall (2012) showed that asymptotic estimates of the residuals variances can be derived which does not involve the parameters and the estimated variances are a function of the empirical distribution function.

The log-logistic density function with shape parameter $\beta > 0$ and scale parameter $\alpha > 0$ is:

$$f(x) = (\beta/\alpha)(x/\alpha)^{\beta-1}/(1+(x/\alpha)^\beta)^2. \qquad (1)$$

The existence of moments and specifically of the mean and variance depends on the shape parameter $\alpha$. The mean is finite for $\beta > 1$ and the variance finite for $\beta > 2$. The distribution function is

$$F(x) = (x/\alpha)^\beta/(1+(x/\alpha)^\beta), \ x > 0. \qquad (2)$$

Some of the applications of the log-logistic distribution is in hydrology (Askhra and Mahdi, 2006), the modelling of incomes in economics, where the distribution is also known as the Fisk distribution (Kleiber and Kotz, 2003) and in survival analysis. It is especially usefull to model accelerated failure times (Klein and Moeschberger, 2003). If $\beta = 1$, it is a generalized Pareto distribution with index two. It is a special case of the Burr distribution (Tadikamalla, 1980).

The cumulative distribution of the two-parameter Weibull distribution is

$$F(x;\alpha,\beta) = 1 - \exp(-(x/\alpha)^\beta), \ x, \alpha, \beta \geq 0, \qquad (3)$$

$\alpha$ a scale parameter and $\beta$ the shape parameter.

## 2 Expressions for the weights used in regression

Expressions for the weights needed in probability plotting type regression, where the order statistics are used will be derived for the two-parameter log-logistic and the Weibull distributions.

### 2.1 The log-logistic distribution

By rewriting (2) as a linear regression form, it follows that:

$$\log(F(x)/(1-F(x))) = \beta \log x - \beta \log \alpha. \qquad (4)$$

Let $x_1, ..., x_n$ be an i.i.d. sample from the log-logistic distribution with corresponding order statistics, $x_{(1)} \leq ... \leq x_{(n)}$. The distribution of $F(x)$ is an uniform distribution on the interval [0,1]. Let $z = \log(F(x)/(1-F(x))$, then it follows that the distribution of $z$ is the logistic distribution with mean 0 and variance $\pi^2/3$, and

$$F(z) = 1/(1+e^{-z}), \quad -\infty < z < \infty.$$

The moments of $z_1 = \log(F(x_{(1)})/(1-F(x_{(1)}))) \leq ... \leq z_n = \log(F(x_{(n)})/(1-F(x_{(n)})))$, would be equal to the moments of the order statistics of a logistic distribution with mean 0 and variance $\pi^2/3$. The derivation of the moments of the logistic distribution is given in the book of Balakrishnan (1991).

Let $\Lambda(x_{(r)}) = \log(F(x_{(r)})/(1-F(x_{(r)})))$, and $\mu_r = E(\Lambda(x_{(r)}))$. Then

$$\Lambda(x_{(r)}) = \beta \log(x_{(r)}) - \beta \log \alpha,$$

$$\Lambda(x_{(r)}) + (\mu_r - \mu_r) = \beta \log(x_{(r)}) - \beta \log \alpha$$

$$\mu_r = \beta \log(x_{(r)}) - \beta \log \alpha + (\mu_r - \Lambda(x_{(r)}))$$

$$\mu_r = \beta \log(x_{(r)}) - \beta \log \alpha + u_r, \qquad (5)$$

where $u_r = \mu_r - \log(F(x_{(r)})/(1-F(x_{(r)})))$, $r=1,...,n$, are the residuals for the regression and the weights are the inverses of the variances of the residuals if only the diagonal elements of the covariance matrix of the residuals are used when performing regression. The density of $F(x_{(r)})$ is a beta distribution with parameters $r$ and $n-r+1$, which does not involve the parameters of the distribution. The moments of $z_r$, $r=1,...n$, can be found by using the moment generating function.

$$M_Z(t) = \frac{1}{B(r, n-r+1)} \int_0^1 e^{t \log(z/(1-z))} z^{r-1} (1-z)^{n-r} dz$$

$$= \frac{1}{B(r, n-r+1)} \int_0^1 z^{t+r-1} (1-z)^{n-r-t} dz$$

$$= \frac{\Gamma(r+t)}{\Gamma(r)} \frac{\Gamma(n-r-t-1)}{\Gamma(n-r+1)}.$$

By noting that $\Gamma'(r+t)/\Gamma(r)|_{t=0} = d\log(\Gamma(r+t))/dt|_{t=0} = \psi(r)$, where $\psi$ denotes the derivative of the log of the gamma function, the digamma function, it follows that

$$\mu_r = E(u_r) = \psi(r) - \psi(n-r+1) \text{ and } Var(u_r) = \psi'(r) + \psi'(n-r+1).$$

Special cases are $r = 1$, $r = n$ and $Var(u_1) = Var(u_n) = \pi^2/6 + \psi'(n)$.

Let $\Sigma$ denote the matrix with $\sigma_{r,s} = Cov(u_r, u_s)$, $r, s = 1,...,n$. It was found that the full covariance matrix has a determinant which is approximately zero when using the asymptotic covariance matrix, and even ridge regression type corrections did not improve estimation result. If it is assumed that the covariance matrix is a diagonal matrix, then it follows that the weight for the $r$-th order statistic is $w_{rr} = 1/(\psi'(r) + \psi'(n-r+1))$ and the sum of squares $\sum_{r=1}^{n} w_{rr} (\mu_r - \beta \log(x_{(r)}) + \beta \log \alpha)^2$, must be minimized with respect to $\alpha$ and $\beta$ to find the weighted least squares estimators of $\alpha$ and $\beta$.

The joint distribution of $v_r = F(x_{(r)})$ and $w_s = F(x_{(s)})$, $r < s$, is

$$f(v_r, w_s) = \frac{v_r^{r-1}(w_s - v_r)^{s-r-1}(1-w_s)^{n-s}}{B(r, s-r)B(s, n-s+1)}, \quad 0 \leq v_s \leq w_r \leq 1. \tag{6}$$

The joint moment generating function of $z_r = \log(v_r/(1-v_r))$ and $z_s = \log(w_s/(1-w_s))$, $r < s$, is

$$M_{r,s:n}(t_1, t_2) = \frac{1}{B(r, s-r)B(s, n-s+1)} \int_0^1 \int_0^{w_s} e^{t_1 z_r + t_2 z_s} v_r^{r-1}(w_s - v_r)^{s-r-1}(1-w_s)^{n-s} dv_r dw_s$$

$$= \frac{1}{B(r, s-r)B(s, n-s+1)} \int_0^1 \int_0^{w_s} v_r^{t_1+r-1}(1-v_r)^{-t_1}(w_s - v_r)^{s-r-1} w_s^{t_2}(1-w_s)^{n-s-t_2} dv_r dw_s.$$

This is exactly the expression when calculating the joint moment generating function of the standard logistic order statistic distribution as given by Balakrishnan (1991). After some algebra the result can be written as

$$M_{r,s:n}(t_1, t_2) = \frac{\Gamma(n+1)}{\Gamma(r)\Gamma(n-s+1)} \sum_{l=0}^{\infty} \frac{(t_1 + l - 1)^{(l)}}{l!} \frac{\Gamma(t_1 + r + l)}{\Gamma(t_1 + s + l)} \frac{\Gamma(t_1 + t_2 + s + l)\Gamma(n - s + 1 - t_2)}{\Gamma(n + t_1 + 1 + l)},$$

$$(t_1 + l - 1)^{(l)} = t_1(t_1 + 1)\ldots(t_1 + l - 1), \quad l \geq 1$$

$$= 1, \quad l = 0.$$

The covariance matrix of the residuals can easily be found by using simulation:

- Simulate say *m* samples of size *n* each, from an uniform distribution.
- Order each sample, say $u_{(1)} \leq ... \leq u_{(n)}$.
- Calculate $z_{(r)} = \log(u_{(r)}/(1-u_{(r)})), r = 1,...,n$, for each sample. This will be equivalent to the ordering of the uniform random numbers, since the transformation is monotone increasing.
- Estimate the covariance matrix from the *m* transformed samples.

This was checked on a few samples and give excellent results. It was found that the covariance matrix is of full rank, but the determinant is approximately zero even for small samples, which can lead to unstable estimation results. For example for *n=25*, the determinant is less than $1.0 \times 10^{-29}$. In the simulation the performance of the estimators will be checked where a diagonal weight matrix is assumed, thus where $W = diag(w_{1,1},...,w_{n,n})$.

$$\mathbf{y}' = [\mu_1,...,\mu_n], \quad X = \begin{pmatrix} 1 & \log(x_{(1)}) \\ \vdots & \vdots \\ 1 & \log(x_{(n)}) \end{pmatrix}.$$

The generalized least squares (GLS) estimate of $\hat{\boldsymbol{\theta}}$, where $\hat{\beta} = \hat{\theta}_2$, $\hat{\alpha} = \exp(-\hat{\theta}_1/\hat{\beta})$ is $\hat{\boldsymbol{\theta}} = (X'WX)^{-1}X'W\mathbf{y}$ and $W = \Sigma^{-1}$. In the simulation the performance of the estimators will be checked where a diagonal matrix is assumed, thus where $W = diag(w_{1,1},...,w_{n,n})$.

It should be noted that this regression is in the inverse regression form and the usual estimators of the variances of the regression coefficients cannot be used, but can easily be approximated by using bootstrap methods.

Van Zyl and Schall (2012) used large sample properties of order statistics, and derived the large sample variance of the residuals. Let $x_p$, where $p$ is defined as $F^{-1}(\xi_p) = p$, and using the asymptotic variance of an order statistic and the large sample variance of a function of a random variable, it follows that:

$$Var(\log(F(x_p)/(1-F(x_p)))) \approx Var(x_p)\left(\frac{d(\log(F(x_p)/(1-F(x_p)))}{dx_p}\right)^2_{|x_p=\xi_p}$$

$$\approx \frac{p(1-p)}{nf^2(\xi_p)}\frac{f^2(\xi_p)}{(p(1-p))^2}$$

$$= 1/(np(1-p)) \ .$$

This approximate variance of $u_r$ can be approximated for the r-th order statistic by:

$$Var(u_r) \approx 1/(n\hat{F}_{(r)}(1-\hat{F}_{(r)})) . \qquad (7)$$

In a similar way it follows that the covariance can be approximated as:

$$Cov(u_r, u_s) = \frac{\hat{F}_r(1-\hat{F}_s)}{nf(x_{(r)})f(x_{(s)})}\frac{f(x_{(r)})f(x_{(s)})}{\hat{F}_r(1-\hat{F}_r)\hat{F}_s(1-\hat{F}_s)}$$

$$= 1/(n\hat{F}_{(s)}(1-\hat{F}_{(r)})). \tag{8}$$

Denote a consistent estimate of the empirical distribution function in the point $x_{(r)}$ by $\hat{F}_r$ and $\hat{\Lambda}(x_{(r)}) = \log(\hat{F}_r/(1-\hat{F}_r))$, for example $\hat{F}_r = r/(n+1)$. The expected value can also be approximated as $\hat{\mu}_r \approx \log(\hat{F}_r/(1-\hat{F}_r))$. Various estimates of $F(x_{(r)})$ for example Bernard's median rank estimator $\hat{F}_r^b = (r-0.3)/(n+0.4)$ (Bernard and Bosi-Levenbach, 1953), can also be used to estimate, $\hat{F}_r^b = \hat{F}(x_{(r)})$.

The focus will be on the log-logistic distribution, but it is easy to show that a similar expression can be derived for the logistic distribution. Consider a logistic distribution with location parameter $\mu$ and scale parameter $\sigma$,

$$F(x) = 1/(1+\exp(-(x-\mu)/\sigma)), \quad -\infty \leq x \leq \infty.$$

It follows that the following regression equation in terms of the order statistics, $x_{(1)} < ... < x_{(n)}$, of a sample of size $n$ from this distribution can be written as:

$$E(\log((1-F(x_{(r)}))/F(x_{(r)}))) = -x_{(r)}/\sigma + \mu/\sigma.$$

Since $Var(\log((1-z_r)/z_r)) = (-1)^2 Var(z_r/(1-z_r))$, it follows that the variance of the residual, $u_r$, corresponding to the r-th order statistic is $Var(u_r) = \psi'(r) + \psi'(n-r+1)$. In a similar way as above it can be shown that the approximate asymptotic variance of $u_r$, is $Var(u_r) = 1/(n\hat{F}_r(1-\hat{F}_r))$ and $Cov(Cov(u_r, u_s)) = 1/(n\hat{F}_s(1-\hat{F}_r))$.

Two estimation methods which are often used to estimate the parameters of the log-logistic distribution are maximum likelihood estimation and probability-weighted moment estimation. For large samples, the maximum likelihood (ML) estimators perform excellent, but it can be shown that better estimators with respect to bias and mean square error (MSE) can be derived. The weighted regression derived in this work outperforms the ML method in small samples, especially with respect to the shape parameter. The large sample properties of the ML estimators were derived by Shoukri *et al.* (1988). The asymptotic biases and variances of the ML estimators, $\hat{\alpha}$ and $\hat{\beta}$ are:

$$bias(\hat{\alpha}) = 1.5\alpha/n\beta^2, \quad bias(\hat{\beta}) = 1.2764\beta/n.$$

The asymptotic covariance matrix of the parameters is, where $\hat{\theta}' = (\hat{\alpha}, \hat{\beta})$:

$$Var(\theta) = \begin{pmatrix} 3\alpha^2/n\beta^2 & 0 \\ 0 & 0.6993\beta^2/n \end{pmatrix}. \tag{9}$$

These expressions can be used to approximate the efficiency of another estimator of $\beta$, say $\beta^*$, as:

$$Eff(\beta^*) = Var(\hat{\beta})/Var(\beta^*). \qquad (10)$$

## 2.2 The Weibull distribution

The cumulative distribution of the two-parameter Weibull distribution is

$$F(x;\alpha,\beta) = 1 - \exp(-(x/\alpha)^\beta), \; x, \alpha, \beta \geq 0,$$

$\alpha$ a scale parameter and $\beta$ the shape parameter. Lieblein 91955) considered the order moments of order statistics of the Weibull distribution. The following equation can be used to estimate the parameters

$$\log(-\log(1-F(x))) = \beta \log(x) - \beta \log(\alpha),$$

The distribution of $F(x)$ is the uniform distribution on the interval [0,1]. Let $z = \log(-\log(1-F(x)))$, then it follows that

$$f(z) = e^z e^{-e^z}, \; -\infty < z < \infty,$$

which is the standard Gumbel distribution for the minimum. The expected value of this distribution is $-\gamma = -0.57721$ (Euler's or Euler–Mascheroni constant). Let $u = z - (-\gamma)$, then it follows that the distribution of the residuals is:

$$f(u) = e^{(u+\gamma)} e^{-e^{(u+\gamma)}}, \quad -\infty < u < \infty, \qquad (11)$$

which is the Gumbel distribution for the minimum with expected value zero and scale parameter 1.

Using similar arguments as for the log-logistic regression the equation used to estimate the parameters with $\mu_r = E(\log(-\log(1 - F(x_{(r)}))))$ is:

$$\mu_r = \beta \log(x_{(r)}) - \beta \log(\alpha) + u_r,$$

$u_r = \mu_r - \log(-\log(1 - F(x_{(r)})))$. In order to find the expressions for the variance of the residuals the following two integrals (Gradshteyn and Ryzhik (1980), p573, p575) are of importance

$$\int_0^\infty e^{-\mu x} \log(x) dx = (-1/\mu)(C + \log(\mu)),$$

and $\int_0^\infty e^{-\mu x} \log^2(x) dx = (1/\mu)(\pi^2/6 + (C + \log(\mu))^2),$

$C$ denotes Euler's constant and is equal to approximately 0.577215.

The expected value of $\log(-\log(1-z_r))$, $z_r = F(x_{(r)})$ and $z_r \sim beta(r, n-r+1)$ is

$$E(\log(-\log(1-z_r))) = (1/B(r,n-r+1)) \int_0^1 \log(-\log(1-z_r)) z_r^{r-1} (1-z_r)^{n-r} dz_r .$$

Let $t = \log(-\log(1-z_r))$, applying the binomial theorem and it follows that

$$E(\log(-\log(1-z_r))) = (1/B(r,n-r+1)) \int_0^\infty \log(t) e^{-t(n-r+1)} (1-e^{-t})^{r-1} dt$$

$$= (1/B(r,n-r+1)) \sum_{j=0}^{r-1} \binom{r-1}{j} (-1)^j \int_0^\infty \log(t) e^{-t(n-r+1+j)} dt, \ r > 1$$

$$= (1/B(r,n-r+1)) \sum_{j=0}^{r-1} \binom{r-1}{j} (-1)^j \left( \frac{-1}{n-r+1+j} \right)$$
$$(C + \log(n-r+1+j)), \ r > 1.$$

In a similar way it can be shown that the second moment for $r > 1$ is

$$E[(\log(-\log(1-z_r)))]^2 = (1/B(r,n-r+1)) \sum_{j=0}^{r-1} \binom{r-1}{j} (-1)^j \left( \frac{1}{n-r+1+j} \right)$$
$$[\pi^2/6 + (C + \log(n-r+1+j))^2], \quad r > 1.$$

For $r = 1$ the expected value is $E(u_1) = -(C + \log(n))$ and $E(u_1^2) = \pi^2/6 + (C + \log(n))^2$,

showing that $Var(u_1) = \pi^2/6$.

In a similar way as for the log-logistic distribution the covariance matrix of the residuals can be found by using simulation:

- Simulate say *m* samples of size *n* each, from an uniform distribution.
- Order each sample, say $u_{(1)} \leq ... \leq u_{(n)}$.
- Calculate $z_{(r)} = \log(-\log(1-u_{(r)}))$, $r=1,...,n$, for each sample. This will be equivalent to the ordering of the uniform random numbers, since the transformation is monotone increasing.
- Estimate the covariance matrix of the transformed sample.

The covariance matrix is of full rank, but for example for *n=10*, the determinant is less than $1.0 \times 10^{-9}$. Thus leading to unstable estimation and it a diagonal weight matrix will be used when performing estimation.

Cohen (1965), White (1963) derived the Cramér-Rao bounds of the ML estimators which are respectively $0.608\beta^2/n$ and $1.109(\alpha/\beta)^2/n$.

# 3 Simulation study results

## 3.1 The accuracy of the approximated residual variances

In the following two plots the true residual variances are plotted together, the estimated values based on 5000 samples, and the approximated residual variances using the asymptotic approximation. The simulated and calculated variances are equal up to at least 0.01 and cannot be distinguished on the charts. The variances of the residuals are functions of the order of the statistics, and not of the parameters, and similar patterns would be found for any parameters of the two distributions considered, the log-logistic and the Weibull distributions.

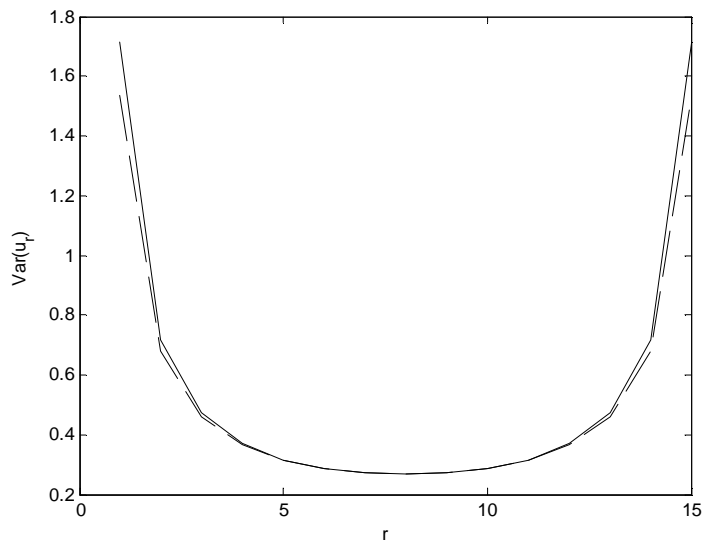

Figure 1. Exact and approximated variances of residuals, log-logistic and n=15 observations.

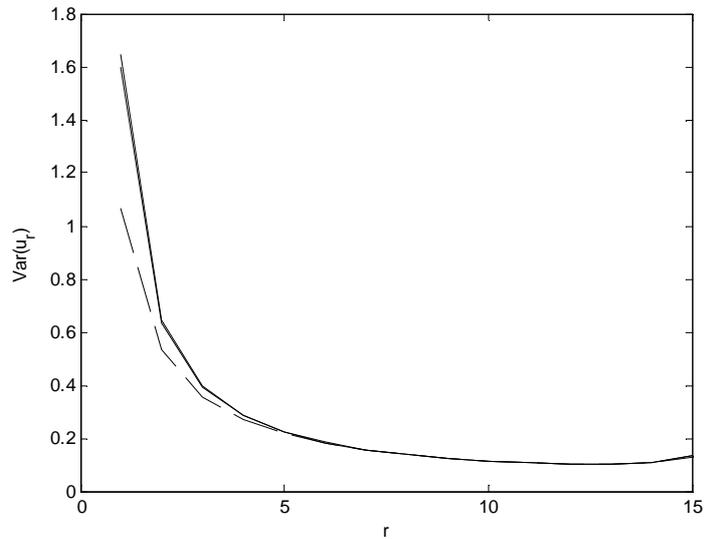

Figure 2. Exact and approximated variances of residuals, Weibull distribution, n=15 observations.

## 3.2 Simulation study concerning the estimation of parameters

In the following two tables estimation was conducted on simulated samples, and the mean square error (MSE) and bias expressed as estimated parameter minus the true parameter are shown. In table 1 samples were generated from the log-logistic distribution. It can be seen that especially in small samples the weighted regression method outperforms ML estimation, especially with respect to bias. The difference in performance was especially with respect to the shape parameter, which might be of more importance in many applications.

|  | $\alpha = 1$ $\beta$ | Weighted LS $\beta$ | Maximum Likelihood $\beta$ | Weighted LS $\alpha$ | Maximum Likelihood $\alpha$ |
|---|---|---|---|---|---|
| **n=15** | **1.0** | -0.0312 (0.0553) | 0.0953 (0.0758) | 0.1060 (0.2913) | 0.1058 (0.2903) |
|  | **1.5** | -0.0380 (0.1292) | 0.1497 (0.1825) | 0.0591 (0.1161) | 0.0587 (0.1158) |
|  | **2.0** | -0.0316 (0.2264) | 0.2213 (0.3224) | 0.0230 (0.0519) | 0.0231 (0.0519) |
|  | **2.5** | -0.0342 (0.3425) | 0.2851 (0.5111) | 0.0103 (0.0323) | 0.0103 (0.0323) |
| **n=25** | **1.0** | -0.0256 (0.0311) | 0.0532 (0.0372) | 0.0684 (0.1525) | 0.0685 (0.1525) |
|  | **1.5** | -0.0496 (0.0718) | 0.0679 (0.0841) | 0.0296 (0.0594) | 0.0296 (0.0595) |
|  | **2.0** | -0.0521 (0.1134) | 0.1024 (0.1321) | 0.0039 (0.0312) | 0.0040 (0.0312) |
|  | **2.5** | -0.0564 (0.1851) | 0.1431 (0.2236) | 0.0157 (0.0194) | 0.0156 (0.0194) |
| **n=50** | **1.0** | -0.0168 (0.0159) | 0.0233 (0.0163) | 0.0403 (0.0673) | 0.0404 (0.0673) |
|  | **1.5** | -0.0237 (0.0330) | 0.0386 (0.0353) | 0.0099 (0.0287) | 0.0099 (0.0287) |
|  | **2.0** | -0.0223 (0.0677) | 0.0594 (0.0734) | 0.0040 (0.0156) | 0.0039 (0.0156) |
|  | **2.5** | -0.0252 (0.0969) | 0.0774 (0.1047) | 0.0063 (0.0098) | 0.0064 (0.0098) |
| **n=100** | **1.0** | -0.0043 (0.0079) | 0.0163 (0.0081) | 0.0166 (0.0318) | 0.0166 (0.0318) |
|  | **1.5** | -0.0144 (0.0169) | 0.0186 (0.0167) | 0.0003 (0.00128) | 0.0002 (0.0128) |
|  | **2.0** | -0.0079 (0.0292) | 0.0315 (0.0296) | 0.0067 (0.0074) | 0.0067 (0.0075) |
|  | **2.5** | -0.0113 (0.0510) | 0.0414 (0.0500) | 0.0012 (0.0052) | 0.0013 (0.0051) |

**Table 1** Results of simulation study based on 1000 samples, log-logistic distribution. Bias and MSE given in brackets.

In table 2 samples were generated from the Weibull distribution. The weighted regression method outperforms ML estimation again in small samples and especially with respect to bias too. The better performance was also in the estimation of the shape parameter with a much smaller bias and also smaller MSE. Even in small samples the ML estimation outperforms the weighted least squares with respect to the estimation of the scale parameter.

|  | $\alpha = 1$ $\beta$ | Weighted LS $\beta$ | Maximum Likelihood $\beta$ | Weighted LS $\alpha$ | Maximum Likelihood $\alpha$ |
|---|---|---|---|---|---|
| n=15 | 0.25 | -0.0066 (0.0031) | 0.0262 (0.0044) | 0.7471 (4.9266) | 0.5473 (3.7747) |
|  | 1.0 | -0.0233 (0.0539) | 0.1017 (0.0761) | 0.0559 (0.0912) | 0.0208 (0.0824) |
|  | 1.5 | -0.0355 (0.1186) | 0.1525 (0.1654) | 0.0125 (0.0341) | -0.0098 (0.0328) |
| n=25 | 0.25 | -0.0066 (0.0019) | 0.0134 (0.0022) | 0.3921 (1.9241) | 0.2751 (1.4957) |
|  | 1.0 | -0.0074 (0.0275) | 0.0723 (0.0353) | 0.0188 (0.0466) | -0.0001 (0.0438) |
|  | 1.5 | -0.0341 (0.0656) | 0.0871 (0.0794) | 0.0079 (0.0188) | -0.0015 (0.0181) |
| n=50 | 0.25 | -0.0028 (0.0009) | 0.0078 (0.0009) | 0.1887 (0.5382) | 0.1405 (0.0461) |
|  | 1.0 | -0.0088 (0.0145) | 0.0301 (0.0149) | 0.0121 (0.0225) | 0.0040 (0.0214) |
|  | 1.5 | -0.0191 (0.0342) | 0.0445 (0.0364) | 0.0105 (0.0112) | 0.0044 (0.0107) |
| n=100 | 0.25 | -0.0020 (0.0005) | 0.0027 (0.0004) | 0.1091 (0.2739) | 0.0788 (0.2437) |
|  | 1.0 | -0.0036 (0.0076) | 0.0169 (0.0068) | 0.0055 (0.0121) | 0.0019 (0.0118) |
|  | 1.5 | -0.0054 (0.0158) | 0.0238 (0.0147) | 0.0003 (0.0050) | -0.0022 (0.0049) |

**Table 2** Results of simulation study based on 1000 samples, Weibull distribution. Bias and MSE given in brackets.

# 4 Conclusions

Exact expressions were derived to estimate the parameters of the log-logistic and Weibull distributions, using weighted regression. The procedure is simple and weights can easily be calculated using simulation. This procedure does not depend on the existence of moments and it was especially effective when estimating the shape parameter.